\documentclass[12pt]{amsart}

\usepackage[latin1]{inputenc}
\usepackage{amsfonts,amsmath,amssymb,amsthm}

\newcommand{\Pp}{\mathbb{P}}

\newcommand{\Nn}{\mathbb{N}}
\newcommand{\Zz}{\mathbb{Z}}
\newcommand{\Qq}{\mathbb{Q}} 
\newcommand{\Ff}{\mathbb{F}} 
\newcommand{\Aa}{\mathbb{A}} 

\newcommand{\DD}{\mathcal{D}} 

\renewcommand{\epsilon}{\varepsilon}
\renewcommand{\le}{\leqslant}
\renewcommand{\ge}{\geqslant}
\renewcommand{\leq}{\leqslant}
\renewcommand{\geq}{\geqslant}

\newcommand{\xbar}{\underline{x}}

\newcommand{\ibar}{\underline{i}}

\def\polsp{S}
\def\spec{\hbox{\rm sp}}
\def\specdiv{\hbox{\rm spdiv}}
\def\card{\hbox{\rm card}}

{\theoremstyle{plain}
\newtheorem{theorem}{Theorem}[section]    
\newtheorem{lemma}[theorem]{Lemma}       
\newtheorem{proposition}[theorem]{Proposition}  
\newtheorem{corollary}[theorem]{Corollary}   
}
{\theoremstyle{remark}
\newtheorem{remark}[theorem]{Remark}   

}

\begin{document}

\title{Indecomposable polynomials and their spectrum}

\author{Arnaud Bodin}
\email{Arnaud.Bodin@math.univ-lille1.fr}

\author{Pierre D\`ebes}
\email{Pierre.Debes@math.univ-lille1.fr}

\author{Salah Najib}
\email{Salah.Najib@math.univ-lille1.fr}

\address{Laboratoire Paul Painlev\'e, Math\'ematiques, Universit\'e 
Lille 1, 59655 Villeneuve d'Ascq Cedex, France}



\date{\today}
 
\begin{abstract}
We address some questions concerning indecomposable 
polynomials and their spectrum. How does the spectrum behave via reduction 
or specialisation, or via a more general ring morphism? Are the indecomposability 
properties equivalent over a field and over its algebraic closure? How many 
polynomials are decomposable over a finite field? 
\end{abstract}

\maketitle


\section{Introduction} 
\label{sec:intro}

Fix an integer $n\geq 2$ and a $n$-tuple of indeterminates $\underline x=(x_1,\ldots,x_n)$. A non-constant polynomial $F(\underline x)\in k[\underline x]$ with coefficients in an algebraically closed field $k$ is said to be {\it indecomposable} in $k[\underline x]$ if it is not of the form $u(H(\underline x))$ with $H(\underline x)\in k[\underline x]$ and $u\in k[t]$ with $\deg(u)\geq 2$. An element $\lambda^\ast\in k$ is called a {\it spectral value} of $F(\underline x)$ if $F(\underline x)-\lambda^\ast$ is reducible in $k[\underline x]$. It is well-known that 
\vskip 2mm

\noindent
(1) {\it $F(\underline x)\in k[\underline x]$ is indecomposable if and only if $F(\underline x) - \lambda$ is irreducible in $\overline{k(\lambda)}[\underline x]$ (where $\lambda$ is an indeterminate)}, 
\vskip 2mm

\noindent
(2) {\it if $F(\underline x)\in k[\underline x]$ is indecomposable, then the subset $\spec(F)\subset k$ of all spectral values of $F(\underline x)$ --- the spectrum of $F(\underline x)$ --- is finite; and in the opposite case, $\spec(F)=k$,}
\vskip 2mm

\noindent
(3) {\it more precisely, if $F(\underline x)\in k[\underline x]$ is indecomposable and for every $\lambda^\ast \in k$, $n(\lambda^\ast)$ is the number of irreducible factors of $F(\underline x) -\lambda^\ast$ in $k[\underline x]$, then we have
$\rho(F) := \sum_{\lambda^\ast\in k} (n(\lambda^\ast)-1)  \leq \deg(F)-1$.
In particular $\hbox{\rm card}(\spec(F))\leq \deg(F)-1$.}
\vskip 2mm

Statement (3), which is known as Stein's inequality, is due to Stein \cite{St} in characteristic $0$ and Lorenzini \cite{Lo} in arbitrary characteristic (but for $2$ variables); see \cite{Na} for the general case.

This paper offers some new results in this context.

In \S \ref{sec:conservation}, given an indecomposable polynomial $F(\underline x)$ with coefficients in an integral domain $A$ and a ring morphism $\sigma:A\rightarrow k$ with $k$ an algebraically closed field,  we investigate the connection between the spectrum of $F(\underline x)$ and that of the polynomial $F^\sigma(\underline x)$ obtained by applying $\sigma$ to the coefficients of $F(\underline x)$. Theorem \ref{thm:main} provides a conclusion {\it \`a la Bertini-Noether}, which, despite its basic nature, does not seem to be available in the literature: under minimal assumptions on $A$, the connection is the expected one generically.
For example if $A=\Zz$, ``spectrum'' and ``reduction modulo a prime $p$'' commute if $p$ is suitably large (depending on $F$). We give other typical applications, notably for a specialization morphism $\sigma$. Related results are given in \cite{BCN}.

For two variables, we can give in \S 3 an in\-de\-com\-po\-sa\-bi\-li\-ty criterion for a reduced polynomial modulo some prime $p$ (theorem \ref{thm:critere petits premiers}) that is more precise than theorem \ref{thm:main}: the condition ``for suitably large $p$'' is replaced by some explicit condition on $F(x,y)$ and $p$, possibly satisfied for small primes. This criterion uses some results on good reduction of curves and covers due to Grothendieck, Fulton et al; we will follow here Zannier's version \cite{Za}. 
Another criterion based on the Newton polygon of a polynomial is given in \cite{CN}.


\S \ref{KvsKbar} is devoted to the connection between the indecomposability properties over a field $K$ and over its algebraic closure $\overline K$. While it was known they are equivalent in many circumstances, for example in characteristic $0$, it remained to handle the inseparable case to obtain a definitive conclusion. That is the purpose of proposition \ref{prop:indecomposabiliteKvsKbar}, which, conjoined with previous works, shows that the only polynomials $F(\underline x)$ indecomposable in $K[\underline x]$ but decomposable in $\overline K[\underline x]$ are $p$-th powers in $\overline K[\underline x]$, where $p>0$ is the characteristic of $K$ (theorem \ref{thm:indecomposabiliteKvsKbar}).

\S \ref{sec:counting} is aimed at counting the number of indecomposable polynomials of a given degree $d$ with coefficients in the finite field $\Ff_q$. We show that most polynomials are indecomposable: the ratio $I_d/N_d$ of indecomposables of  degree $d$ tends to $1$ (as $d\rightarrow \infty$ or as 
$q\rightarrow \infty$), and we give some estimate for the error term $1 - I_d/N_d$ 
The constants involved in our estimates are explicit 
(as in \cite{vzG1} for irreducible polynomials). 
For simplicity we mostly restrict to polynomials in two variables as calculations become more intricate when $n>2$. We also consider the one variable situation (for which the definition of indecomposability is slightly different, see \S \ref{ssec:one variable 1}) with the restriction that $q$ and $d$ are relatively prime. The cases ($n>2$) and 
($n=1$ with $(q,d)\not=1$) are considered by J. von zur Gathen in a parallel work \cite{vzG2} \cite{vzG3}. 
\medskip

\noindent
{\bf Acknowledgments.} We wish to thank J. von zur Gathen for interesting discussions in Lille
and in Bonn.


\section{Spectrum and morphisms} \label{sec:conservation}

{\bf Notation:} 
If $\sigma: A\rightarrow B$ is a ring morphism, we denote  the image of elements $a\in A$ by $a^\sigma$. For $P(\underline x)\in A[\underline x]$, we denote the polynomial obtained by applying $\sigma$ to the coefficients of $P$ by  $P^\sigma(\underline x)$. If $V\subset \Aa^n_{A}$ is the Zariski closed subset associated with a family of polynomials $P_i(\underline x) \in A[\underline x]$, we denote by $V^\sigma$ the Zariski closed subset of $\Aa^n_{B}$ associated with the family of polynomials $P_i^\sigma(\underline x) \in B[\underline x]$. 

If $S\subset A$ is a multiplicative subset such that all elements from $S^\sigma$ are invertible in $B$, we still denote by $\sigma$ the natural extension $S^{-1}A\rightarrow B$ of the original morphism $\sigma$. 

\smallskip

Fix an integrally closed ring $A$, with a perfect fraction field $K$.

An effective divisor $D=\sum_{i=1}^r n_i a_i$ of $\overline K$ is said to be {\it $K$-rational} if the coefficients  of the polynomial $P(T)=\prod_{i=1}^r (T-a_i)^{n_i}$ are in $K$\footnote{which, under our hypothesis ``$K$ perfect'', amounts to the invariance of $P(T)$, or of $D$, under $\hbox{\rm Gal}(\overline K/K)$.}. A morphism $\sigma: A\rightarrow k$ in an algebraically closed field $k$ is then said to be {\it defined at $D$} if the coefficients of $P(T)$ have a common denominator $d\in A$ such that $d^\sigma$ is non-zero in $k$ \footnote{which, under our hypothesis ``$A$ integrally closed'', amounts to saying the elements $a_i$ themselves have a common denominator $d\in A$ (that is, $da_i$ integral over $A$, $i=1,\ldots,r$) such that $d^\sigma\not= 0$.}. In this case we denote by $P^\sigma(T) \in k[T]$ the image polynomial of $P(T)$ by the morphism $\sigma$ (extended to the fraction field of $A$ with denominators a power of $d$) and by $D^\sigma$ the effective divisor
  of $k$ whose support is the set of roots of $P^\sigma(T)$ and coefficients are the corresponding multiplicities.

\subsection{Statement} \label{ssec: statements}
For more precision, we use the {\it spectral divisor} rather than the spectrum: it is the divisor $\specdiv(F) = \sum_{\lambda^\ast\in k} (n(\lambda^\ast)-1) \hskip 2pt \lambda^\ast$ of the affine line $\Aa^1(k)$. Its support is the spectrum 
of $F$ and  Stein's inequality 
rewrites: $\deg(\specdiv(F)) \leq \deg(F)-1$.

\begin{theorem} \label{thm:main} Let $F(\underline x)\in A[\underline x]$ be indecomposable in $\overline K[\underline x]$. 
Then there exists a non-zero element 
$h_F\in A$ such that the following holds. For every morphism $\sigma:A\rightarrow k$ in an algebraically closed field $k$, if $h_F^\sigma\not=0$, then $F^\sigma(\underline x)$ 
is indecomposable in $k[\underline x]$, the morphism $\sigma:A\rightarrow k$ is defined at the divisor $\specdiv(F)$ and we have $\specdiv(F^\sigma)=(\specdiv(F))^\sigma$; in particular $\rho(F^\sigma)=\rho(F)$ and $\spec(F^\sigma)=(\spec(F))^\sigma$.





\end{theorem}

The first stage of the proof will produce the spectrum as a Zariski closed subset of the affine line $\Aa^1_A$ over the ring $A$. Specifically the following can be drawn from the proof:  there is a proper\footnote{that is, distinct from the whole surrounding space (here the affine line $\Aa^1_A$ over the ring $A$); equivalently, there exists a non-zero polynomial in the associated ideal.}  Zariski closed subset $V_F\subset \Aa^1_A$ such that for every morphism $\sigma:A\rightarrow k$ as above,
\vskip 1mm

\noindent
{\rm (*)} {\it the polynomial $F^\sigma(\underline x)$, if it is of degree $d$, is indecomposable 
in $k[\underline x]$ if and only if the Zariski closed subset $V_F^\sigma \subset \Aa^1_k$ is proper, 
and in this case, we have $\spec(F^\sigma) = V_F^\sigma(k)$}.
\vskip 1mm

\noindent
When applied to the inclusion morphism $A\rightarrow \overline K$, theorem \ref{thm:main} yields that the spectrum of $F(\underline x)$ is equal to the Zariski closed subset $V_F(\overline K)$. In particular, it is $K$-rational. The same is true for the spectral divisor of $F(\underline x)$ as $n(\lambda^\tau) = n(\lambda)$ for each $\lambda \in \overline K$ and each $\tau \in \hbox{\rm Gal}
(\overline K/K)$.
 
%

\subsection{Typical applications}

\subsubsection{Situation 1} For $A=\Zz$, then $h_F\in \Zz$, $h_F\not=0$.  
Theorem \ref{thm:main}, applied with
$\sigma:\Zz\rightarrow \overline \Ff_p$ the reduction morphism modulo a prime number $p$ yields:
 
\smallskip

\noindent
{\it for all suitably large $p$, the reduced polynomial $\overline F(\underline x)$
modulo $p$ is indecomposable in $\overline{\Ff_p}[\underline x]$ and its spectral
divisor is obtained by reducing that of $F(\underline x)$, that is: 
 $\specdiv(\overline F) = \overline{\specdiv(F)}$.}

\subsubsection{Situation 2} Take $A=k[\underline t]$ with $k$ an algebraically closed field and $\underline t = (t_1,\ldots,t_r)$ some indeterminates. Denote in this situation by  $F(\underline t,
\underline x)$ the polynomial $F(\underline x)$ of the general statement. Theorem \ref{thm:main}, applied with $\sigma$ the specialisation morphism $k[\underline t]\rightarrow k$ that maps
$\underline t = (t_1,\ldots,t_r)$ on an $r$-tuple $\underline t^\ast=
(t_1^\ast,\ldots,t_r^\ast)\in k^r$ yields:
\smallskip

\noindent
{\it for all $\underline t^\ast$ off a proper Zariski closed subset of $k^r$, 
the specialized polynomial  $F(\underline t^\ast,\underline x)$
is indecomposable in $k[\underline x]$ and its spectral divisor  is obtained
by specializing that of $F(\underline t,\underline x)$.}

\subsubsection{Situation 3} $F(\underline x)$ is the generic polynomial in $n$ variables and of degree $d$. Take for $A$ the ring $\Zz[a_{\ibar}]$ generated by the indeterminates $a_{\ibar}$ 
corresponding to the coefficients of $F(\underline x)$; the multi-index ${\ibar}=(i_1,\ldots,i_n)$ ranges over the set $I_{n,d}$ of all $n$-tuples of integers $\geq 0$ such that $i_1+\cdots+i_n \leq d$. 

Classically the polynomial $F(\underline x)$ is
irreducible in $\overline {\Qq(a_{\ibar})}[\underline x]$, hence it is  indecom\-po\-sa\-ble.
Theorem \ref{thm:main}, applied with $\sigma:A\rightarrow k$ a specialization morphism of the $a_{\ibar}$, yields that all polynomials $f(\underline x)\in k[\underline x]$ of degree $d$ are indecomposable but possibly those from the proper Zariski closed subset corresponding to
the equation $h_F = 0$ (with $h_F$ viewed in $k[a_{\ibar}]$).

For polynomials $f(\underline x)$ outside the closed subset $h_F = 0$, the spectrum of $f$
is obtained by specializing the generic spectrum. However we have:

\begin{proposition} \label{prop:spectre non vide}
For $d>2$ or $n>2$, the generic spectrum is empty. For $d=2$, it contains a single element, given by
\vskip -1mm

$$a_{00} - \frac{a_{02} a_{10}^2 + a_{20} a_{01}^2 - a_{01}a_{10}a_{11}}{4 a_{02}a_{20}-a_{11}^2}$$
\end{proposition}

For $d>2$ or $n>2$, polynomials with a non-empty spectrum lie in the Zariski closed subset $h_F = 0$. 

\begin{proof} Assume that the generic spectrum is not empty. If $k$ is an algebraically closed field and ${\mathcal R}_{n,d}$ (resp. ${\mathcal P}_{n,d}|_{a_{\underline 0}=0}$) denotes the set of polynomials $P(\underline x)\in k[\underline x]$ of degree $\leq d$ that are reducible in $k[\underline x]$ (resp. 
whose constant term is zero), the correspondence $P(\underline x)\rightarrow P(\underline x)-P(\underline 0)$ induces an algebraic morphism ${\mathcal R}_{n,d}\rightarrow{\mathcal P}_{n,d}|_{a_{{\underline 0}}=0}$ which is generically surjective (that is, surjective above a non-empty Zariski open subset of  ${\mathcal P}_{n,d}|_{a_{{\underline 0}}=0}$). It follows that ${\mathcal R}_{n,d}$ is of codimension $\leq 1$ in
the space ${\mathcal P}_{n,d}$ of all polynomials in $k[\underline x]$ of degree $\leq d$. 
This observation provides the desired conclusion in the case $n=2$ and $d>2$: indeed we have
 $\hbox{\rm codim}_{{\mathcal P}_{2,d}}({\mathcal R}_{2,d}) = d-1$ \cite[theorem 2]{vzG1}. 

For $d=2$, the equation ``$(ux+ay+b)(vx+cy+d)=F(\underline x)$ modulo the constant term'' with unknowns $u, a, b, v, c, d$ is readily solved: reduce to the case $a_{20}=u=v=1$, find the unique solution for the $4$-tuple $(a,b,c,d)$ and compute $bd$; the generic spectral value is then $a_{00}-bd$.

Finally assume that for $d\geq 2$ and $n>2$, there exists a generic spectral value $\lambda \in \overline K$ (with $K= \Qq(a_{\ibar})$). Let $F(\underline x) -\lambda = Q(\underline x) R(\underline x)$ be a  non trivial factorization  in $\overline K[\underline x]$. Specializing $x_3,\ldots,x_n$ to $0$ gives a non trivial factorization in $\overline K[x_1,x_2]$ of the generic polynomial of degree $d$ in $2$ variables. From the first part of the proof, we have $d=2$. Furthermore, the above case provides the necessary value of $\lambda$. Now specializing $x_2$ and $x_4,\ldots,x_n$ to $0$ leads to a different value. Whence a contradiction.
\end{proof}

\subsection{Proof of theorem \ref{thm:main}}

\subsubsection{1st stage: elimination theory} \label{etape1} This stage is aimed at showing proposition \ref{prop:bertini-noether generalized} below, which generalizes the Bertini-Noether theorem \cite[prop.9.4.3]{FrJa}.
It is proved in the general situation
\vskip 2mm 

\noindent
(Hyp) {\it a polynomial ${\mathcal F}(\underline{\lambda}, \underline x)\in A[\underline{\lambda}, \underline x]$ is irreducible in $\overline {K(\underline{\lambda})}[\underline x]$, where $\underline{\lambda} = (\lambda_1,\ldots,\lambda_s)$ is an $s$-tuple of indeterminates ($s\geq 0$).} 
\vskip 2mm

\noindent
We will use it in the special case  ${\mathcal F}(\underline{\lambda}, \underline x)= F(\underline x)-\lambda$. The hypotheses ``$A$ integrally closed''  and ``$K$ perfect'' are not necessary for this stage.

As in situation 3, consider  some indeterminates $(a_{\ibar})_{{\ibar}\in I_{n,d}}$ corresponding
to the coefficients of a polynomial of degree $d$ in $n$ variables. A polynomial with coefficients in a ring $R$ corresponds 
to a morphism $\phi: \Zz[a_{\ibar}]\rightarrow R$; denote by
$F(a_{\ibar}^\phi)(\underline x) \in A[\underline x]$ the corresponding polynomial. 
Let $\varphi_{\underline \lambda}: \Zz[a_{\ibar}]\rightarrow A[\underline \lambda]$ be the morphism corresponding to the polynomial from statement (Hyp): ${\mathcal F}(\underline \lambda, \underline x) = F(a_{\ibar}^{\varphi_{\underline \lambda}})(\underline x)$. 

From Noether's theorem \cite[\S 3.1 theorem 32]{Sc}, there exist finitely many universal homogeneous forms ${\mathcal N}_{h}(a_{\ibar})$ ($1\leq h \leq D=D(n,d)$) in the $a_{\ibar}$ and with coefficients in $\Zz$ such that: 
\smallskip

\noindent
(4) {\it for every morphism $\phi:\Zz[a_{\ibar}]\rightarrow k$ 
in an algebraically closed field $k$, the polynomial 
$F(a_{\ibar}^\phi)(\underline x)$,  if it is of degree $d$, is reducible in 
$k[\underline x]$ if and only if ${\mathcal N}_{h}(a_{\ibar}^\phi)=0$ 
for $h=1,\ldots,D$.}

\smallskip
For $\phi$ taken to be the morphism $\varphi_{\underline \lambda}: \Zz[a_{\ibar}]\rightarrow A[\underline \lambda] \subset \overline{K(\underline \lambda)}$, the elements ${\mathcal N}_{h}(a_{\ibar}^{\varphi_{\underline \lambda}})\in A[\underline \lambda]$ are polynomials $N_{h}(\underline \lambda)$. Let $V_{\mathcal F}\subset \Aa_A^s$ be the Zariski closed subset corresponding to the ideal they generate; it is a proper closed subset. Indeed, as ${\mathcal F}(\underline{\lambda}, \underline x)$ is irreducible  in $\overline {K(\underline{\lambda})} [\underline x]$, from (4), at least one of the polynomials $N_{h}(\underline \lambda)$, say $N_{h_0}(\underline \lambda)$, is non-zero. Denote by $a_{\mathcal F}\in A$ the product of a non-zero coefficient of $N_{h_0}(\underline \lambda)$ and the non-zero coefficient of some monomial of ${\mathcal F}(\underline{\lambda}, \underline x)$ of degree $d$ in $\underline x$.

If $R$ is an integral domain and $\Sigma:A[\underline{\lambda}]\rightarrow R$ a morphism, then
(4), with $\phi$ taken to be $\Sigma \circ \varphi_{\underline \lambda}: \Zz[a_{\ibar}]\rightarrow R \hookrightarrow \kappa$ and $\kappa=\overline{{\rm Frac}(R)}$, yields
that the polynomial ${\mathcal F}^\Sigma \in R[\underline x]$, if of degree $d$, is irreducible in $\kappa[\underline x]$ if and only if at least one of the elements $N_{h}^\Sigma \in R$
is non-zero (note that ${\mathcal F}^\Sigma = F(a_{\ibar}^{\Sigma\varphi_{\underline \lambda}})(\underline x)$ and ${\mathcal N}_{h}(a_{\ibar}^{\Sigma \varphi_{\underline \lambda}})= {\mathcal N}_{h}(a_{\ibar}^{\varphi_{\underline \lambda}})^\Sigma$), or, equivalently, if the corresponding Zariski closed subset of ${\rm Spec}(R)$ 
is proper. 

Let $\sigma:A\rightarrow k$ be a morphism with $k$ algebraically closed. Apply the above
first to the morphism $\sigma \circ \varphi_{\underline \lambda}: \Zz[a_{\ibar}]\rightarrow k[\underline \lambda]$ and then, for $\underline{\lambda}^\ast \in k^s$, to the morphism $s_{\underline{\lambda}^\ast}\circ  \sigma \circ \varphi_{\underline{\lambda}}:\Zz[a_{\ibar}] \rightarrow k$ obtained
by composing $\sigma \circ \varphi_{\underline{\lambda}}$ with the specialization morphism $s_{\underline{\lambda}^\ast}: k[\underline \lambda] \rightarrow k$ to   $\underline \lambda^\ast$. Conclude:

\begin{proposition}[Bertini-Noether generalized]  \label{prop:bertini-noether generalized}     \hskip 3mm

{\rm (a)} The polynomial ${\mathcal F}^\sigma(\underline{\lambda}, \underline x)$, if it is of degree $d$ in $\underline x$, is irreducible in $\overline{k(\underline{\lambda})}[\underline x]$ if and only if the Zariski closed subset $V_{\mathcal F}^\sigma \subset \Aa_{k}^s$ is proper. All these conditions are satisfied if $a_{\mathcal F}^\sigma$ in non-zero in $k$.

{\rm (b)} If the polynomial ${\mathcal F}^\sigma(\underline{\lambda}^\ast, \underline x)$ is of degree $d$, then it is reducible in $k[\underline x]$ if and only if $\underline{\lambda}^\ast$ is in the set ${V}_{\mathcal F}^\sigma(k)$. 
\end{proposition}

\subsubsection{2nd stage: implications for the spectrum of $F(\underline x)$}  \label{etape2}
We return to the situation where ${\mathcal F}(\lambda,\underline x) = F(\underline x)-\lambda$. Denote the Zariski closed subset $V_{\mathcal F}$ from \S \ref{etape1} by $V_F$; it is a Zariski closed subset of the affine line $\Aa^1_A$. The preceding conclusions, conjoined with the connection, recalled in \S \ref{sec:intro}, between indecomposability of  $F(\underline x)$ and irreducibility of $F(\underline x)-\lambda$, yield statement (*) from \S \ref{ssec: statements}. 

Denote by $s_F(\lambda)$ the g.c.d.\ of the polynomials $N_{h}(\lambda)$ in the ring $K[\lambda]$. Write it as  $s_F(\lambda) = \polsp_F(\lambda)/c_1$ with
$\polsp_F(\lambda) \in  A[\lambda]$ and $c_1 \in A$ non-zero.
The polynomial $\polsp_F(\lambda)$ is non-zero and its distinct roots in $\overline K$, say $\lambda_1,\ldots,\lambda_s$, which are the common roots in $\overline K$ of the polynomials $N_h(\lambda)$, are the spectral values of $F(\underline x)$ (note that $F(\underline x)-\lambda^\ast$ is of degree $d$ for all $\lambda^\ast\in \overline K$). Thus we have $S_F(\lambda) = c_2 \prod_{i=1}^s (\lambda-\lambda_i)^{n_i} \in A[\lambda]$ for some exponents $n_i>0$ and $c_2\in A$, $c_2\not=0$. It follows that the set $\spec(F)=\{\lambda_1,\ldots,\lambda_s\}$ is $K$-rational. As already noted, the same is then true for the spectral divisor $\specdiv(F)$. 

\subsubsection{3rd stage: invariance of the spectrum of $F$ via morphisms }  \label{etape3}
Fix a morphism $\sigma:A\rightarrow k$ with $k$ algebraically closed. Denote by $a_F$ the element $a_{\mathcal F}$ from \S \ref{etape1} for ${\mathcal F}=F(\underline x)-\lambda$. If $a_F^\sigma\not= 0$, $F^\sigma(\underline x)$ is of degree $d$ and indecomposable in $k[\underline x]$. 
Furthermore, its spectral values are the roots in $k$ of the g.c.d. of the polynomials $N_{h}^\sigma(\lambda)$. 

Note that the element $c_2$ above is a common denominator of $\lambda_1,\ldots,\lambda_s$; if $c_2^\sigma\not= 0$, the morphism $\sigma:A\rightarrow k$ is defined at $\specdiv(F)$. 

\begin{lemma} \label{lemme:pgcd et reduction} There exists $c_3\in A$, $c_3\not=0$ such that, if $a_F^\sigma c_1^\sigma c_2^\sigma c_3^\sigma\not= 0$, the polynomial $S_F^\sigma(\lambda)\in k[\lambda]$ equals (up to some non-zero multiplicative constant in $k$) the g.c.d. in $k[\lambda]$ of polynomials $N_{h}^{\sigma}(\lambda)$ ($1\leq h \leq D$). In particular $\spec(F^\sigma)=(\spec(F))^\sigma$.
\end{lemma}

\begin{proof}
The problem is whether the g.c.d. commutes with $\sigma$. The Euclidean algorithm
provides the g.c.d. as the last non-zero remainder. To reach our goal, it suffices to guarantee
that for each division $a=bq+r$ in $K[\lambda]$ involved in the algorithm, the identity 
$a^\sigma=b^\sigma q^\sigma+r^\sigma$, with $\sigma$ suitably extended, be the 
division of  $a^\sigma$ by $b^\sigma$ in $k[\lambda]$. For this, write $a$, $b$, $q$ as $r$
in the form $n(\lambda)/m$ with $n(\lambda)\in A[\lambda]$ and $m\in A$, consider
the product $\beta$ of denominators $m$ of $a$, $b$, $q$ and $r$ with the coefficients 
of highest degree monomials in the numerators $n(\lambda)$ of $b$ and $r$ and request
that $\beta^\sigma\not= 0$. Multiplying all elements $\beta$ for all divisions 
leading to the g.c.d. of two, then of all polynomials
in question, leads to a non-zero element $c_3 \in A$ which satisfies the desired
statement.
\end{proof}

\begin{remark} Morphisms and g.c.d. do not commute in general: for example 
$\gcd(\lambda, \lambda+a)$ is $1$ generically, but equals $\lambda$ if $a= 0$.
\end{remark}


\subsubsection{4th stage: invariance of $\specdiv(F)$ via morphisms}  \label{etape4} It remains to extend the conclusion ``$\spec(F^\sigma)=(\spec(F))^\sigma$'' to the spectral divisor $\specdiv(F)$. We will show how to guarantee that,  {\it via} the morphism $\sigma$, the spectral values remain distinct and the associated decompositions of $F(\underline x)-\lambda$ have the same numbers of distinct irreducible factors\footnote{The argument will also show the degrees of these irreducible factors, say $Q_{\lambda,j}$, remain the same and thus so does the quantity $\min_{\lambda\in \spec(F)} (\sum_j  \deg(Q_{\lambda,j}) - 1$ which replaces $\deg(F) - 1$ in Lorenzini's refined version \cite{Lo} of Stein's inequality.}.

Consider  the discriminant of the polynomial $\prod_{i=1}^s (\lambda-\lambda_i)$; it is a non-zero element of $K$. Write it as $c_4/c_5$ with $c_4,c_5\in A$, non-zero. If $c_4^\sigma c_5^\sigma\not=0$, 
the polynomials $S_F(\lambda)$ and $S_F^\sigma(\lambda)$ have the same number of distinct roots, whence $\hbox{\rm card}(\spec(F^\sigma))=\hbox{\rm card}((\spec(F))^\sigma)=\hbox{\rm card}(\spec(F))$.

For $i=1,\ldots,s$, let $F(\underline x) -\lambda_i= \prod_{j=1}^{n(\lambda_i)} Q_{ij}(\underline x)^{k_{ij}}$ be a factorization (into distinct irreducible polynomials) in $\overline K[\underline x]$. Let $E/K$ be a finite Galois extension that contains the finite set ${\mathcal C}$ of all coefficients involved in all above factorizations, $c_6$ be a non-zero element of $A$ such that $c_6 c$ is integral over $A$ for all
$c\in {\mathcal C}$ and $c_7$ be the discriminant of a basis of $E$ over $K$ the elements of which are integral over $A$. Denote by $B$ the fraction ring of $A$ with denominator a power of $c_6c_7$ and by $B^\prime_E$ the integral closure of $B$ in $E$. The ring $B^\prime_E$ is a free $B$-module of rank $[E:K]$. Assume that $c_6^\sigma c_7^\sigma\not=0$. The morphism $\sigma:A\rightarrow k$ extends to a morphism $B\rightarrow k$, and, as $k$ is algebraically closed, this morphism $\sigma:B\rightarrow k$ can in turn be extended to a morphism $\widetilde \sigma:B^\prime_E\rightarrow k$. 

The polynomials $Q_{ij}(\underline x)$ are in the ring $B^\prime_E[\underline x]$ and are absolutely irreducible. The (classical) Bertini-Noether theorem provides a non-zero element $\beta\in B^\prime_E$ such that, if $\beta^{\widetilde \sigma}\not=0$, then each of the polynomials $Q_{ij}^{\widetilde \sigma}(\underline x)$ is absolutely irreducible. Therefore the decomposition $F^\sigma(\underline x) -\lambda_i^{\widetilde \sigma}= \prod_{j=1}^{n(\lambda_i)} Q_{ij}^{\widetilde \sigma}(\underline x)$ obtained from the preceding one by applying $\widetilde \sigma$, is the factorization of $F^\sigma(\underline x) -\lambda_i^{\widetilde \sigma}$ into  irreducible polynomials in $k[\underline x]$. 

It remains to assure that for $i$ fixed, the polynomials $Q_{ij}^{\widetilde \sigma}(\underline x)$ are different, even up to non-zero multiplicative constants. For any two (distinct) polynomials $Q_{ij}(\underline x)$, $Q_{ij^\prime}(\underline x)$, the matrix with rows the tuples of coefficients of the two polynomials has a $2\times 2$-block with a non-zero determinant. Denote the product of all such determinants for all possible couples $(Q_{ij}(\underline x),Q_{ij^\prime}(\underline x))$ by $\delta$; 
it is a non-zero element of $B^\prime_E$. Denote then by $\nu$ the norm of $\beta\delta$ relative to the extension $E/K$. As $A$ is integrally closed, so is $B$ and $\nu\in B$. 
Write it as  $\nu=c_8/(c_6c_7)^\gamma$ with $c_8\in A$ and $\gamma\in \Nn$. Condition
 $c_6^\sigma c_7^\sigma c_8^\sigma\not=0$ implies $\beta_F^{\widetilde \sigma}
 \delta_F^{\widetilde \sigma}\not=0$. Theorem \ref{thm:main} 
 is finally established for $h_F= a_F \prod_{i=1}^8 c_i$.
 
 \begin{remark} The same proof, with the polynomial ${\mathcal F}(\underline {\lambda}, \underline x)$ from \S \ref{etape1} of the form $F(\underline x)-\lambda G(\underline x)$ with $F(\underline x)$, $G(\underline x)\in A[\underline x]$ and $\deg G \le \deg F$, leads to the more general form of theorem \ref{thm:main} for which
 indecomposable polynomials are replaced by indecomposable rational functions (in this case, ``indecomposable'' means not of the form $u(H(\underline x))$ with $H(\underline x)$ and $u(t)$ rational functions and $\deg(u)\geq 2$ \footnote{the degree of a rational function is the maximum of the degrees of its numerator and denominator.}). A spectral value of a rational function $F(\underline x)/G(\underline x)$ is an element $\lambda$ such that the polynomial $F(\underline x)-\lambda G(\underline x)$ is reducible. Statements (1), (2) and (3)
from \S  \ref{sec:intro} remain true, except that the bound in Stein's inequality should be replaced by $(\deg(F))^2 -1$ \cite{Bo} \cite{Lo}.
More generally one can take ${\mathcal F}(\underline {\lambda}, \underline x)$ of the form $F(\underline x)-\lambda_1 G_1(\underline x) - \cdots -\lambda_s G_s(\underline x)$ with $F(\underline x), G_1(\underline x),\ldots,G_s(\underline x) \in A[\underline x]$ and handle other situations studied in the
literature. In this context,  some effective results are given in \cite{BCN}. 
\end{remark}

\section{An indecomposability criterion modulo $p$}\label{sec:critere}

In this section $n=2$, $A$ is a Dedekind domain and its fraction field $K$ is assumed to be of characteristic $0$. Fix also a non-zero prime ideal $\mathfrak{p}$ of $A$ and assume its residue field $k=A/\mathfrak{p}$ is of characteristic $p>0$. Denote by $\widetilde x$ the image of an element $x$ by the reduction morphism $A\rightarrow k$. The situation ``$A=\Zz$ and $\mathfrak{p}=p\Zz$'' is typical. 

Let $F(x,y)\in A[x,y] $ be an indecomposable polynomial in $\overline K[x,y]$ of degree 
$d\geq 1$, monic in $y$. 

Here is our strategy to guarantee indecomposablity of $F(x,y)$ 
modulo $\mathfrak{p}$. Pick $\lambda^\ast \in A\setminus \spec(F)$ (using Stein's theorem, this can be done with $\lambda^\ast$ not too big). Thus $F(x,y)-\lambda^\ast$ is irreducible
in $\overline K[x,y]$. It follows from the classical Bertini-Noether theorem that if ``$\mathfrak{p}$ is
big enough'', then the reduced polynomial $F(x,y)-\lambda^\ast$ modulo $\mathfrak{p}$ is absolutely
irreducible. Therefore $F(x,y)$ is indecomposable modulo $\mathfrak{p}$ (as there is at least one
non spectral value). However the constants involved in the condition ``$\mathfrak{p}$ big enough'' are
too big for a practical algorithmic use. We will follow an alternate approach, based
on good reduction criteria for covers, and more precisely Zannier's criterion \cite{Za}.

Consider the discriminant with respect to $y$ of $F(x,y)-\lambda$: 
\vskip -2mm

$$\Delta_F(x,\lambda)=\hbox{\rm disc}_y(F(x,y)-\lambda)$$ 

\noindent
Denote then the product of all distinct irreducible factors of $\Delta_F(x,\lambda)$ in $K(\lambda)[x]$
by $\Delta_F^{\rm red}(x,\lambda)$; more precisely, $\Delta_F^{\rm red}(x,\lambda)$ is defined by the following formula, which is also algorithmically more practical:
\vskip -1mm

$$\Delta_F^{\rm red}(x,\lambda) = c(\lambda) \hskip 2mm \frac{\Delta_F(x,\lambda)}{\hbox{\rm gcd}(\Delta_F(x,\lambda),(\Delta_F)^\prime_x(x,\lambda))}$$ 

\noindent
where the g.c.d.\ is calculated in the ring $K(\lambda)[x]$ (using the Euclidean algorithm for example) and $c(\lambda)\in K(\lambda)$ is the rational function, defined up to some invertible element in $A$, 
that makes $\Delta_F^{\rm red}(x,\lambda)$ a primitive polynomial in $A[\lambda][x]$. Consider next 
the polynomial:
$$\Delta_F(\lambda) = \hbox{\rm disc}_x(\Delta^{\rm red}_F(x,\lambda)).$$
\noindent
We have $\Delta_F(\lambda)\in A[\lambda]$ and $\Delta_F(\lambda)\not=0$. Finally let $\Delta_0(\lambda)\in A[\lambda]$
be the coefficient of the highest monomial in $\Delta_F(x,\lambda)$ (viewed in $A[\lambda][x]$).

\begin{theorem} \label{thm:critere petits premiers} Assume, in addition to $F(x,y)$ being indecomposable in $\overline K[x,y]$, that the reduced polynomial $\widetilde \Delta_0(\lambda)\widetilde\Delta_F(\lambda)$ is non-zero in $k[\lambda]$ and that $p>\deg_Y(F)$. Then $\widetilde{F}(x,y)$ is indecomposable in $\overline k[x,y]$.
\end{theorem}

The assumption $p>\deg_Y(F)$ can be replaced by the weaker condition that $p$ does not divide the order of the Galois group of $F(x,y)-\lambda$, viewed as a polynomial in $\overline{K(\lambda)}(x)$ 
(see footnote 8). 

\begin{remark}
Theorem \ref{thm:critere petits premiers} can be combined with preceding results. If $V_F\subset \Aa^1_A$ is the Zariski closed subset from \S \ref{ssec: statements}, then, under the above hypotheses, the reduced Zariski closed subset $\widetilde{V_F}\subset \Aa^1_{\overline k}$ is proper and its points are the spectral values of $\widetilde F$: $\spec(\widetilde F)=\widetilde{V_F}(\overline k)$. However the assumptions on $p$ and $F(x,y)$ may not be sufficient to guarantee the extra conclusions
$\spec(\widetilde F )=\widetilde{\spec(F)}$ and $\specdiv(\widetilde{F})=\widetilde{\specdiv(F)}$ from theorem \ref{thm:main} (which may not even be well-defined).
\end{remark}


\begin{proof}[Proof of theorem \ref{thm:critere petits premiers}]The prime ideal $\mathfrak{p} \subset A$ determines a discrete valuation $v$ of $K$ whose valuation ring is the localized ring $A_{\mathfrak{p}}$; the fraction field of $A_{\mathfrak{p}}$ and its residue field remain equal to $K$ and $k$ respectively. Hypotheses and conclusions from theorem \ref{thm:critere petits premiers} are unchanged if $A$ is replaced by $A_{\mathfrak{p}}$. The valued field $(K,v)$ can then also be replaced by any finite extension of the completion $K_v$ and $A$ by the new valuation ring; the discrete valuation $v$ uniquely extends, the residue field is replaced by some (finite) extension of $k$, the indecomposability properties of $F(x,y)$ over $K$ or over $K_v$ are equivalent. 

Thus we may and will assume that $(K,v)$ is a complete discretely valued field, that $A$ is its valuation ring (which is integrally closed) and that the field $K$ and the residue field $k$ contain as many  (finitely many) algebraic elements over the original fields as necessary.

The polynomial $\Delta_F(x,\lambda) $ is in $A[x,\lambda]$ and its factorization into irreducible polynomials in $K(\lambda)[x]$ can be written
\vskip -2mm

$$\Delta_F(x,\lambda) = \delta_0(\lambda) \prod_{i=1}^s \Delta_i(x,\lambda)^{\alpha_i}$$
\vskip -1mm

\noindent
where the polynomials $\Delta_i(x,\lambda)$ are in $A[x,\lambda]$, irreducible in $K(\lambda)[x]$, pairwise distinct (even up to some constant in $K$) and are primitive in $A[\lambda][x]$, where $\delta_0(\lambda) \in A[\lambda]$ and where the $\alpha_i$ are positive integers. Then, up to some invertible element in $A$, we have
\vskip -2mm

$$\Delta_F^{\rm red}(x,\lambda) =  \prod_{i=1}^s \Delta_i(x,\lambda)$$

\noindent
Also note that the polynomial $\Delta_0(\lambda)$ is a multiple in $A[\lambda]$ of the product of $\delta_0(\lambda)$ with the highest monomial coefficients $\delta_1(\lambda),\ldots,\delta_s(\lambda)$ of the polynomials $\Delta_1(x,\lambda),\ldots,\Delta_s(x,\lambda)$ (viewed in $A(\lambda)[x]$).

Pick next $\widetilde\lambda^\ast \in k$ such that $\widetilde \Delta_0(\widetilde \lambda^\ast )\widetilde\Delta_F(\widetilde \lambda^\ast )\not=0$ in $k$, then lift it to some element $\lambda^\ast\in A$ such that $\lambda^\ast\notin \spec(F)$. This is possible in view of the preliminary remark. 

The set of roots of $\Delta_F(x,\lambda^\ast)$ contains the set of finite\footnote{{\it i.e.}, distinct from the point at infinity.} branch points of the cover of $\Pp^1_x$ \footnote{The subscript ``$x$''  indicates that the cover is induced by the correspondence $(x,y)\rightarrow x$. In fact the problem is symmetric in the variables $x$ and $y$ which can be switched in our statement.} determined by the (absolutely irreducible) polynomial $F(x,y)-\lambda^\ast$. 
The preliminary remark makes it possible to assume that these roots are in $K$. Furthermore as $\widetilde \delta_i (\widetilde\lambda^\ast )\not=0$, we have $\delta_i(\lambda^\ast)\in A\setminus \mathfrak{p}$, $i=1,\ldots,s$; therefore these roots are integral over $A$ and so are in 
$A$.

As $\Delta_F(\lambda^\ast)\not=0$, the roots of $\Delta_F^{\rm red}(x,\lambda^\ast)$ in $\overline K$ are  distinct and as $\delta_0(\lambda^\ast)\not=0$, they are the roots of $\Delta_F(x,\lambda^\ast)$. As $\widetilde \Delta_0(\widetilde \lambda^\ast )\not=0$, $\widetilde \Delta_F(x,\widetilde \lambda^\ast)$ is not the zero polynomial. As $\widetilde \Delta_F(\widetilde \lambda^\ast )\not=0$, the roots of $\widetilde\Delta_F^{\rm red}(x,\widetilde\lambda^\ast)$, which are those of the polynomial $\widetilde\Delta_F(x,\widetilde \lambda^\ast )$, are distinct. Thus we obtain that the distinct roots of the polynomial $\Delta_F(x,\lambda^\ast)$, and {\it a fortiori} the branch points  of the cover considered above, 
have distinct reductions modulo the ideal $\mathfrak{p}$.

It follows from standard results on good reduction of covers, and more precisely here, from the main theorem of \cite{Za} that, under the assumption $p>\deg_Y(F)$ \footnote{It suffices to assume that $p$ does not divide the order of the Galois group of $F(x,y)-\lambda^\ast$, which divides the order of the Galois group of $F(x,y)-\lambda$, which itself divides $(\deg_Y(F))!$ .}, $\widetilde {F}(x,y)-\widetilde \lambda^\ast$ is absolutely irreducible. Hence $\widetilde F(x,y)$ is indecomposable in $\overline k[x,y]$. 
\end{proof}


\section{Indecomposability over $K$ versus $\overline K$} \label{KvsKbar}

\subsection{Statements \hbox{\rm (for $n\geq 2$ variables)}} The indecomposability property which we recalled the definition of  in \S \ref{sec:intro} over an algebraically closed field  can in fact be defined over an arbitrary field: just require that the polynomials $u(t)$ and $H(\underline x)$ involved have their coefficients in the field in question. The results below identify the only cases where the property is not the same over some field $K$ and over some extension $E$. The following result handles the case that $E/K$ is purely inseparable, which was missing in the literature.\smallskip

\begin{proposition} \label{prop:indecomposabiliteKvsKbar} Let $E/K$ be a purely inseparable 
algebraic field extension of characteristic $p>0$ and $F(\underline x)\in K[\underline x]$. Assume $F(\underline x)$ is 
not of the form $b\hskip 2pt G(\underline x)^p + c$ with $G(\underline x)\in E[\underline x]$ 
and $b,c\in K$. Then $F(\underline x)$ is indecomposable in $K[\underline x]$ 
if and only if it is indecomposable in $E[\underline x]$.
\end{proposition}

If $E=\overline K$, the assumption on $F(\underline x)$ rewrites to merely say that $F(\underline x)$ is not a $p$-th power in $\overline K[\underline x]$, which in turn is equivalent to at least one exponent in $F(\underline x)$ not being a multiple of $p$. Clearly this assumption cannot be removed: for example, if $\alpha \in \overline K\setminus K$ but $\alpha^p=a\in K$ then $x^p+ay^p$ is indecomposable in $K[\underline x]$ but decomposable in $\overline K[\underline x]$.

In \cite[proposition 1]{ArPe}, Arzhantsev and Petravchuk show the equivalence from proposition \ref{prop:indecomposabiliteKvsKbar} without any assumption on $F(\underline x)$, but in the case of a separable extension $E/K$ (possibly of positive transcendence degree).
As any extension is a purely inseparable algebraic extension of some separable extension, 
conjoining their result with ours yields that, under the assumption on $F(\underline x)$ from proposition \ref{prop:indecomposabiliteKvsKbar}, the equivalence holds for an arbitrary extension $E/K$. We can be more precise.

\begin{theorem} \label{thm:indecomposabiliteKvsKbar} Let $E/K$ be a field extension and $F(\underline x)\in K[\underline x]$ be a non-constant polynomial. Then the following are equivalent:
\vskip 2mm

\noindent
{\rm (i)} $F(\underline x)$ is indecomposable in $K[\underline x]$ but decomposable in $E[\underline x]$.
\vskip 2mm

\noindent
{\rm (ii)} {\rm (a)} $K$  is of characteristic $p>0$ and  $E/K$ is inseparable, 
\vskip 1mm

\hskip 3mm {\rm (b)} $F(\underline x) = b\hskip 2pt G(\underline x)^p + c$ for some $G(\underline x)\in E[\underline x]$ and $b,c\in K$, and
\vskip 1mm

\hskip 3mm {\rm (c)}  $G(\underline x)^p$ is indecomposable in $K[\underline x]$.

\end{theorem}

Condition (ii) (c) implies that $G(\underline x)$ is not of the form $u(H(\underline x))$ with $u\in E[t]$, $H(\underline x)\in E[\underline x]$, $\deg(u)\geq 2$ and both $u(t)^p\in K[t]$ and $H(\underline x)^p\in K[\underline x]$. But there are other possible polynomials that should be excluded whose description is more intricate.

\subsection{Proofs} \label{ssec:proofs KvsKbar}

\begin{proof}[Proof of proposition \ref{prop:indecomposabiliteKvsKbar}] The converse part is obvious. For the direct part, assume $F(\underline x)$ is decomposable in $E[\underline x]$. Then  it is decomposable over some finite extension of $K$ contained in $E$, which admits a finite system of generators $\alpha_1,\ldots,\alpha_s$ with irreducible polynomial over $K$ of the form $x^{p^n}-a$
with $a\in K$. The multiplicativity of the degree and of the separable degree imply that the extensions $K(\alpha_1,\ldots,\alpha_{j+1})/ K(\alpha_1,\ldots,\alpha_{j})$ are purely inseparable, $j=1,\ldots,s-1$.
By induction one reduces to the case $s=1$, and then a new induction reduces to the case $E=K(\alpha)$ with $\alpha^p = a\in K\setminus K^p$.

Assume $F(\underline x) = h(G(\underline x))$ with $h(t)\in K(\alpha)[t]$ such that 
$\deg(h)\geq 2$ and $G(\underline x)\in K(\alpha)[\underline x]$. We deduce

$$F(\underline x)^p  = {}^ph(G(\underline x)^p)$$

\noindent
where, if $h(t)=\sum_{i=0}^{\deg(h)} h_it^i$, we set ${}^ph(t) = \sum_{i=0}^{\deg(h)} h_i^pt^i$. 
As ${}^ph(t) \in K[t]$ and $G(\underline x)^p\in K[\underline x]$ (since $y^p\in K$ for all $y\in K(\alpha)$), this shows that the field $K(F(\underline x), G(\underline x)^p)$ is of transcendence degree $1$ over $K$. From Gordan's theorem \cite[\S 1.2, th.3]{Sc}, there exists $\theta(\underline x)\in K(\underline x)$ such that

$$K(F(\underline x), G(\underline x)^p)=K(\theta(\underline x))$$

\noindent
Furthermore from \cite[\S 1.2, th.4]{Sc}, one may assume that $\theta(\underline x)\in K[\underline x]$. 
Thus we have

$$\left\{\begin{matrix}
F(\underline x) = u(\theta(\underline x)) \ \hbox{\rm with}\ u(t)\in K(t) \hfill \\
G(\underline x)^p = v(\theta(\underline x)) \ \hbox{\rm with}\ v(t)\in K(t) \hfill \\
\end{matrix}
\right.$$

\noindent
As $F(\underline x)$ and $G(\underline x)^p$ are polynomials, $u(t), v(t)$ are necessarily in $K[t]$. 
It follows from the indecomposability of $F(\underline x)$ over
$K$ that $\deg(u) = 1$, which gives $G(\underline x)^p = 
w(F(\underline x))$ for some polynomial $w\in K[t]$. But then we obtain $G(\underline x)^p = w\circ h(G(\underline x))$, which, since $G(\underline x)$
is non constant, amounts to $T^p=w\circ h(T)$ where $T$ is an indeterminate.
As $\deg(h)\geq 2$ and $p$ is a prime, we have $\deg(w)=1$ and $\deg(h)=p$, which
gives $F(\underline x) = b\hskip 2pt G(\underline x)^p + c$ for some $b,c\in K$. 

Note that because of the inductive process, conclusion ``$b,c\in K$'' should really be that 
$b,c$ are in the first subfield of the initial reduction. But $F(\underline x)$ being in $K[\underline x]$ then implies that $b\gamma^p \in K$ for some non-zero $\gamma\in E$ and $b\hskip 2pt G(\underline 0)^p+c\in K$. Up to changing $G(\underline x)$ to $\gamma^{-1}G(\underline x) - \gamma^{-1}G(\underline 0)$,
one can indeed conclude that $b,c\in K$ in the general situation.
\end{proof}

\begin{proof}[Proof of theorem \ref{thm:indecomposabiliteKvsKbar}]
{\rm (i) $\Rightarrow$ (ii)}: If $K_s/K$ is the maximal separable extension contained in $E$, then, from the Arzhantsev-Petravchuk result, $F(\underline x)$ is indecomposable in $K_s[\underline x]$. In particular $E\not=K_s$, which gives (ii) (a). Proposition \ref{prop:indecomposabiliteKvsKbar} then provides condition (ii) (b) except that $b$ and $c$ are {\it a priori} in $K_s$, but using again the final note of the proof of Proposition \ref{prop:indecomposabiliteKvsKbar}, one can indeed choose $b,c\in K$. Condition (ii) (c) then readily follows from (ii) (b) and the indecomposability of $F(\underline x)$ in $K[\underline x]$. The other implication {\rm (ii) $\Rightarrow$ (i)} is clear.
\end{proof}

\subsection{One variable} \label{ssec:one variable 1}
In proposition \ref{prop:indecomposabiliteKvsKbar}, $F(\underline x)$ is a polynomial in two variables or more. In one variable, the indecomposability definition 
should be modified (for otherwise it is trivial): a polynomial $F(x)\in k[x]$ is said to be indecomposable in $k[x]$ if it is not of the form $u(H(x))$ with $H(x)\in k[x]$ and $u\in k[t]$ with $\deg(u)\geq 2$ \emph{and} $\deg(H)\geq 2$.

\begin{proposition} Proposition \ref{prop:indecomposabiliteKvsKbar} holds for one variable polynomials.
\end{proposition}

\begin{proof}
The same proof can be used as for proposition \ref{prop:indecomposabiliteKvsKbar}. It leads to

$$\left\{\begin{matrix}
F(x) = u(\theta(x)) \ \hbox{\rm with}\ u(t)\in K[t] \hfill \\
G(x)^p = v(\theta(x)) \ \hbox{\rm with}\ v(t)\in K[t] \hfill \\
\end{matrix}
\right.$$

\noindent
But from the indecomposability of $F(\underline x)$ over $K$, we now deduce that
$\deg(u) = 1$ or $\deg(\theta) = 1$.

The case $\deg(u) = 1$ is handled as before. In the other case, we deduce 
from $\deg(\theta) = 1$ that $K(F(x),G(x)^p)=K(x)$, which implies that 
$K(\alpha)(h(G(x)),G(x)^p) = K(\alpha)(x)$ and so that 

$$K(\alpha)(x)\subset K(\alpha)(G(x))$$

\noindent
which forces $\deg(G)=1$ and contradicts the decomposability assumption in one variable made at the beginning of the proof.
\end{proof}


\section{Counting indecomposable polynomials over finite fields} \label{sec:counting}
For each integer $d\geq 1$, denote the number of polynomials in $\Ff_q[\underline x]$ 
($\underline x=(x_1,\ldots,x_n)$)
of degree $d$ by $N_d$. 
We have $$\left\{
\begin{matrix}
& N_d =  \left( q^{\binom{n+d-1}{n-1}} -1 \right) \cdot q^{\binom{n+d-1}{n}}\ \hbox{\rm (for general $n$) }\hfill\\
& N_d = q^{\frac 12 (d+1)(d+2)}(1-q^{-d-1})\ \hbox{\rm (for $n=2$)} \hfill \\
& N_d= (q-1)q^d \ \hbox{\rm (for $n=1$)} \hfill \\
\end{matrix}
\right.$$ 

\noindent
Denote the number of those polynomials which are indecomposable 
(resp.~decomposable) by $I_d$ (resp. $D_d$). We have $N_d=I_d+D_d$.

We will study separately the case of $n\geq 2$ variables (\S \ref{ssec:main result} - \S \ref{ssec:proof of (c)}) and the case $n=1$ (\S \ref{ssec:one variable (counting)}).

\subsection{Main result} \label{ssec:main result} From \S \ref{ssec:main result} to \S \ref{ssec:proof of (c)}, we assume $n\geq 2$. 

\begin{theorem} \label{thm:counting,n=2} 
\noindent
{\rm (a)} $I_d/N_d$ tends to $1$ in the two situations where $d \to \infty$ with $q$ fixed,
and where $q \to \infty$ with $d$ fixed.
\vskip 2mm

\noindent
{\rm (b)} If $d$ is a product of at most $2$ prime numbers $p\leq p^\prime$, then

\begin{itemize}
  \item $d=p$ and $D_d = q^{d}(q^n-1)$, or
  \item $d = p^2$  and $ D_{d} = q^{p-1}N_p+(q^{d}-q^{2p-1})(q^n-1)$, or
  \item $d = p p^\prime$ with $p<p^\prime$ and 
  
  $ D_{d} = q^{p-1}N_{p^\prime}+ q^{p^\prime-1}N_p+(q^{d}-2q^{p+p^\prime-1})(q^n-1)$.
   \end{itemize}
\vskip 2mm

\noindent   
 {\rm (c)}  Assume $n=2$.  If $d$ is the product of at least $3$ prime numbers, then
   
   $$ \left| \frac{D_d}{N_d} - \alpha_{d} \right| \le \alpha_{d} \hskip 2pt \beta_{d}\hskip 5mm  \hbox{\rm where}\ \displaystyle \left\{\begin{matrix} \displaystyle
   \alpha_{d} = \frac{q^{\ell-1+\frac 12 (\frac{d}{\ell}+1)(\frac{d}{\ell}+2)}}{q^{\frac 12 (d+1)(d+2)}}\hfill\\
\displaystyle   \beta_{d} =  \frac{d}{q^{\frac{d}{\ell}}}\hfill\\
   \end{matrix}\right.$$
   
   \noindent
and $\ell>1$ is the first (hence prime) divisor of $d$.
 \end{theorem}

\subsection{An induction formula} \label{ssec:induction formula}

Let $K$ be an arbitrary field. Let $F = u\circ H$ be a decomposition of $F \in K[\xbar]$  with $u \in K[t]$, $\deg u \ge 2$, and $H \in K[\xbar]$. We say that $F = u \circ H$ is a \emph{normalized decomposition} if $H$ is indecomposable, monic ({\it i.e.} the coefficient of the leading term of a chosen order is $1$) and its constant term equals zero. Given a decomposition $F=u\circ H$, there exists an associated 
normalized decomposition $F = u^\prime \circ H^\prime$. The following lemma shows it is unique.

\begin{lemma}
\label{prop:uni}
Let $ F=u \circ H = u^\prime\circ H^\prime$ be two normalized decompositions of $F\in K[\xbar]$. 
Then $ u=u^\prime$ and $H=H^\prime$.
\end{lemma}

\begin{proof}
It follows from $u(H)-u^\prime(H^\prime)=0$ that $H$ and $H^\prime$ are algebraically dependent
over $K$.  By Gordan's theorem \cite[\S 1.2, theorems 3 and 4]{Sc} (already used in \S \ref{ssec:proofs KvsKbar}),
there exists a polynomial $\theta(\xbar) \in K[\xbar]$ such that $K[\theta] = K[H,H^\prime]$.  
That is, there exist $v,v^\prime \in K[t]$ such that $H=v(\theta)$ and $H^\prime=v^\prime(\theta)$.
As the two decompositions of $F$ are normalized, $H$ and $H^\prime$ are
indecomposable, so $\deg v=\deg v^\prime=1$, and so using the other normalization
conditions, we obtain $H = H^\prime$. Finally it follows from $u(H)=u^\prime(H)$ that $u=u^\prime$.
\end{proof}

\begin{corollary}[induction formula] \label{cor:induction formula} With notation as in \S \ref{ssec:main result}, we have 
$$ I_d = N_d - \sum_{d^\prime | d\hskip 2pt, \hskip 2pt d^\prime< d} 
q^{\frac{d}{d^\prime}-1}\times I_{d^\prime}$$
\end{corollary}

\begin{proof}
Let $d^\prime\geq 1$ be a divisor or $d$. There are $(q-1)q^{d/d^\prime}$ polynomials $u\in \Ff_q[t]$ of degree $d/d^\prime$ and $I_{d^\prime}/q(q-1)$ normalized indecomposable polynomials $H\in \Ff_q[\xbar]$ of degree $d^\prime$. The formula follows as from lemma \ref{prop:uni}, every 
polynomial $F$ counted by $D_d$ can be uniquely written $F= u \circ H$ with $u$ and $H$ as above 
for some integer $d^\prime$ such that $d^\prime | d\hskip 2pt, \hskip 2pt d^\prime< d$.
\end{proof}

Conjoined with $I_1=N_1= q(q^n-1)$ this formula provides an algorithm to compute
$I_d$ and $D_d$, which is convenient for small $d$.

\subsection{Proof of theorem \ref{thm:counting,n=2}  \hbox{\rm (a)} and  \hbox{\rm (b)}} \label{ssec:proof of (a) and (b)}
The formulas in (b) straightforwardly follow from corollary \ref{cor:induction formula}. If $d=p$ is a prime number, we have
$D_p = q^{p-1}  I_1 = q^{p-1} N_1 = q^{p} (q^n-1)$.
If $d=p^2$ then
  \begin{align*}
    D_d &= q^{p-1}  I_p + q^{p^2-1}  I_1 \\
    &= q^{p-1}  (N_p-q^{p}(q^n-1)) + q^{p^2}(q^n-1). \\
  \end{align*}

\vskip -6mm
\noindent
Computations are similar for $d=p p^\prime$.
To prove (a) we write
 
 $$N_d-I_d = D_d = \sum_{d^\prime | d\hskip 2pt, \hskip 2pt d^\prime< d} q^{d/d^\prime }I_{d'} \le \sum_{d^\prime |d \hskip 2pt, \hskip 2pt d^\prime< d} q^{d/d^\prime}N_{d^\prime}$$
The sum has at most $d$ terms and each is $\leq q^{d}N_{d/2}$, whence

$$1-\frac{I_d}{N_d} \le d \hskip 2pt q^{d} \hskip 2pt \frac{N_{d/2}}{N_d}$$

\noindent
and the announced result as the right-hand side term tends to $0$ in the two situations considered in the statement
of theorem \ref{thm:counting,n=2} (a). $\square$

\subsection{Proof of theorem \ref{thm:counting,n=2} \hbox{\rm (c)}} 
\label{ssec:proof of (c)} In this subsection we assume that $n=2$ and that $d$ has at least three prime divisors.

\subsubsection{A technical lemma} 

\begin{lemma}
\label{lem:bd}
Let $b(d) = \frac12 (d+1)(d+2)$.
Let $\ell>1$ be the first divisor of $d$ and $\ell^\prime>\ell$ be the second divisor of $d$.
Let $\lambda \ge \ell^\prime$ be a divisor of $d$ and $\ell''>1$ be the first divisor of $d/\ell$.
Then we have
\begin{enumerate}
  \item \label{it:l1} $b(d/\ell^\prime)+\ell^\prime\ge b(d/\lambda)+\lambda$.
  \item \label{it:l2} $b(d/\ell) + \ell - d/\ell \ge b(d/\ell^\prime)  + \ell^\prime$.
  \item \label{it:l3} $b(d/\ell) + 1 - d/\ell \ge b(d/\ell\ell'')+\ell''$.
\end{enumerate}
\end{lemma}

\begin{proof} 
(1) We have
$$b(d/\ell^\prime)+\ell^\prime - b(d/\lambda)-\lambda = \frac 12 
\left(\frac{d}{\ell^\prime}-\frac{d}{\lambda} 
\right)\left(\frac{d}{\ell^\prime}+\frac{d}{\lambda} + 3 - 
2\frac{\ell^\prime\lambda}{d} \right) \ge 0$$
because $\displaystyle d/\ell^\prime -d/\lambda \ge 0$ and $ 
\displaystyle \frac{d}{\ell^\prime}+\frac{d}{\lambda} + 3 - 2\frac{\ell^\prime\lambda}{d} \ge 
\frac{d}{\ell^\prime}+4-2\ell^\prime \ge 0$ as $d$ has at least $3$ prime divisors.
\vskip 2mm

(2) We have $\ell \ell^\prime \le d$ so $\displaystyle \ell^\prime-\ell \le \frac d \ell$.
Moreover we have $\displaystyle \frac d {\ell^\prime} \le \frac d \ell -2$ and for all $d 
\ge 6$ we have $b(d/\ell^\prime) \le b(d/\ell - 2)$.
Hence
$$
 b(d/\ell) - b(d/\ell^\prime)+ \ell - \ell^\prime - 
\frac{d}{\ell} \ge b(d/\ell)- b(d/\ell - 2)
-2d/\ell = 1
$$

(3) If we set $\delta = d/\ell$ then
$$b(\delta)+1-\delta - b(\delta/\ell'')-\ell'' = \frac 12 
\left(\delta-\frac{\delta}{\ell''} 
\right)\left(\delta+\frac{\delta}{\ell''} -2 \right)
+ \frac12 \left(3\delta-5 \frac{\delta}{\ell''}-2\ell''+2 \right)$$

\noindent
Now $\delta-\frac{\delta}{\ell''} \ge 0$, $\delta+\frac{\delta}{\ell''} -2 
\ge 0$
and as $\delta$ has at least $2$ prime divisors,  then $u(\ell'') = 3\delta-5 
\frac{\delta}{\ell''}-2\ell''+2 \ge u(2) = \frac \delta 2  - 2 \ge 0$.
\end{proof}

\subsubsection{An upper bound for $D_d$} \label{ssec:upperbound} Using the notations of Lemma \ref{lem:bd}, we have
 \begin{align*}
D_{d}    &= q^{\ell-1} I_{d/\ell} + \sum_{\lambda | d, \lambda >\ell} q^{\lambda-1} I_{d/\lambda}  \qquad \text{(corollary \ref{cor:induction formula})}\\
    &\le q^{\ell-1} N_{d/\ell} + \sum_{\lambda | d, \lambda >\ell} q^{\lambda-1} N_{d/\lambda}  \\
    &\le q^{b(d/\ell)+\ell-1}\left(1-\frac 1 {q^{\frac d \ell+1}}\right)+  \sum_{\lambda | d, \lambda >\ell}q^{\lambda-1}q^{b(d/\lambda)} \qquad \text{(explicit formula for $N_{d/\lambda}$)} \\
    &\le q^{b(d/\ell)+\ell-1}\left(1-\frac 1 {q^{\frac d \ell+1}}\right)+  (d-1) q^{b(d/\ell^\prime)+\ell'-1} \qquad \text{(lemma \ref{lem:bd} (\ref{it:l1}))} \\
    &\le q^{b(d/\ell)+\ell-1}\left(1-\frac 1 {q^{\frac d \ell+1}}\right) \left( 1+  \frac{d}{q^{b(d/\ell) -b(d/\ell^\prime) + \ell - \ell^\prime}} \right) \qquad \text{(because } \frac{d-1}{1-q^{-\frac d \ell -1}} \le d) \\
    &\le q^{b(d/\ell)+\ell-1}\left(1-\frac 1 {q^{\frac d \ell+1}}\right) \left( 1+  \frac{d}{q^{\frac{d}{\ell}}} \right) \qquad \text{(lemma \ref{lem:bd} (\ref{it:l2})) } \\
  \end{align*}
  
\subsubsection{A lower bound for $D_d$}
Start from $D_{d} \ge  q^{\ell-1} I_{d/\ell} $.
Then use \S \ref{ssec:upperbound} right above (or the formulas already proved from theorem \ref{thm:counting,n=2} (b)) to bound $I_{d/\ell}= N_{d/\ell} - D_{d/\ell}$ from below. We obtain 

  \begin{align*}
    D_{d} &\ge q^{\ell-1}\times \left( q^{b(\frac{d}{\ell})}\left(1-\frac 1 {q^{\frac d \ell +1}}\right)   -  
   q^{b(\frac{d}{\ell\ell''})+\ell''-1} \left(1-\frac 1 {q^{\frac{d}{\ell\ell''}+1}}\right) \left( 1 + \frac{d/\ell}{q^{\frac{d}{\ell\ell''}}}  \right)  \right) \\
   &\ge  q^{\ell-1} \left(1-\frac 1 {q^{\frac d \ell +1}}\right) \left( q^{b(\frac{d}{\ell}) } -  
    2 q^{b(\frac{d}{\ell\ell''})+\ell''-1} \right) \qquad (\text{because } \frac{d/{\ell}}{q^{\frac{d}{\ell\ell''}}} \le 1) \\
        &= q^{\ell-1}\left(1-\frac 1 {q^{\frac d \ell +1}}\right) q^{b(d/\ell)} \left(1- \frac{2}{q^{ b(d/\ell) -b(d/\ell\ell'') +1-\ell''  }}  \right) \\
    &\ge  q^{b_{\frac{d}{\ell}}+\ell-1}\left(1-\frac 1 {q^{\frac d \ell+1}}\right)\left( 1 -  \frac{2}{q^{\frac{d}{\ell}}} \right) \qquad \text{(lemma  \ref{lem:bd}  (\ref{it:l3}))} \\
  \end{align*}

\subsubsection{Final estimate for $D_d/N_d$.}

The upper and lower bounds for $D_d$ yield the following inequalities
$$\frac{q^{b(\frac{d}{\ell})+\ell-1}}{q^{b(d)}} \times \frac{1-q^{-\frac d \ell-1}}{1-q^{-d-1}} \times \left( 1- \frac{2}{q^{\frac d \ell}} \right) \le \frac{D_d}{N_d} \le \frac{q^{b(\frac{d}{\ell})+\ell-1}}{q^{b(d)}} \times \frac{1-q^{-\frac d \ell-1}}{1-q^{-d-1}} \times \left( 1+ \frac{d}{q^{\frac d \ell}} \right)$$
\noindent
which are a little more precise than the announced statement. $\square$

\subsection{One variable} \label{ssec:one variable (counting)}

Here we assume $n=1$. For polynomials in one variable, we use the definition of indecomposability given in \S \ref{ssec:one variable 1}.

\subsubsection{Main result}
\begin{theorem}\label{th:one}
Assume $q$ and $d$ are relatively prime.

\noindent
{\rm (a)} If $d$ is a product of at most $2$ prime numbers $p\leq p^\prime$, then
\begin{itemize}
  \item $d=p$ and $D_d=0$, or
  \item $d=p^2$ and $D_{d} = \frac{q-1}{q}\hskip 1pt q^{2p}$, or
  \item $d=pp^\prime$ with $p<p^\prime$ and
$$  2\frac{q-1}{q}q^{p+p^\prime} - q^5 \le D_{d} \le 2\frac{q-1}{q}q^{p+p^\prime}$$
\end{itemize}

\noindent
{\rm (b)}  Assume $d$ is the product of at least $3$ prime numbers. Let $\ell>1$ be the first 
divisor of $d$ and $\ell^\prime >\ell$ be its second divisor. Then 
we have 
\vskip -2mm

$$\frac{d}{2\ell}\frac{1}{q^{\frac{d}{\ell}-\frac{d}{\ell^2}-\ell+1}} \le \frac{D_d}{N_d} - \alpha_{d} \le \frac{d-2}{2q^{\ell + \frac{d}{\ell}-\ell^\prime- \frac{d}{\ell^\prime}}}\hskip 5mm  \hbox{where} \hskip 5mm \alpha_{d} = \frac{2}{q^{d-\ell-\frac{d}{\ell}+1}}$$
\end{theorem}

As a consequence we have that $I_d/N_d$ tends to $1$ in the two situations where $d \to \infty$ with $q$ fixed, and $q \to \infty$ with $d$ fixed.

Theorem \ref{th:one} fails if the assumption $(q,d)=1$ is removed. For example for $q=2$ and $d$ even one can compute that $\displaystyle D_d/N_d \sim 3.2^{-d/2}$ while $\alpha_d = 4. 2^{-d/2}$ in this case. 

From now on we assume $q$ and $d$ are relatively prime.
The rest of the paper is devoted to the proof of theorem \ref{th:one}. Our strategy is similar to the one used for $n\geq 2$. We view the set $\DD_d$ of all decomposable polynomials $f(x)\in \Ff_q[x]$ of degree $d$ as the union of smaller sets which we will estimate.  More specifically we write 
$$ \DD_d = \bigcup_{\lambda | d\hskip 2pt, \hskip 2pt \ell \le \lambda \le d/\ell} \DD_{\lambda,d/\lambda}$$ 
where $\DD_{\lambda,d/\lambda}\subset \DD_d$ is the subset of all 
$f(x)$ which admit a decomposition $f = u \circ v$ with $u,v\in \Ff_q[x]$,  $\deg u = \lambda \ge 2$, $\deg v = d/\lambda \ge 2$,  $v$ monic and of constant term equal to $0$.
A difference with the case $n\geq 2$ is that we do \emph{not} have a partition.

\subsubsection{1st stage: upper bounds} \label{ssec:first stage} (Assumption $(q,d)=1$ is not used in this paragraph). For every divisor $\lambda\geq 1$ of $d$, denote the cardinality of $\DD_{\lambda,\frac{d}{\lambda}}$ by $D_{\lambda,\frac{d}{\lambda}}$. 
We have 
\vskip -2mm

$$D_{\lambda,\frac{d}{\lambda}} \leq N_\lambda \hskip 2pt \frac{N_{d/\lambda}}{q(q-1)}= \frac{q-1}{q} \hskip 2pt q^{\lambda+\frac{d}{\lambda}}$$

\noindent
If $\ell>1$ is the first divisor of $d$ and $\ell^\prime>\ell$ the second divisor, we have
\vskip -2mm

$$  D_d \le \sum_{\lambda | d\hskip 2pt, \hskip 2pt \ell \le \lambda \le d/\ell} {D_{\lambda,\frac d \ell}} \le \frac{q-1}{q}   \sum_{\lambda | d\hskip 2pt, \hskip 2pt\ell \le \lambda \le d/\ell} q^{\lambda + \frac{d}{\lambda}}$$

\noindent
The idea is that the main contribution comes from $D_{\ell,\frac{d}{\ell}}$ and $D_{\frac{d}{\ell},\ell}$. 

If $d$ is the product of exactly $2$ prime numbers $\ell$ and $d/\ell$, then these are the only contributions and we have the desired upper bound.
Otherwise we write $\lambda + \frac{d}{\lambda} \le \ell^\prime + \frac{d}{\ell^\prime}$ to bound the extra terms 
and obtain
$$D_d \le  \frac{q-1}{q}  q^{\ell + \frac{d}{\ell}}\left( 2 + \frac{d-2}{q^{\ell + \frac{d}{\ell}-\ell^\prime - \frac{d}{\ell^\prime}}} \right)$$
which yields all announced upper bounds in theorem \ref{th:one}. We also deduce this practical bound: $\displaystyle D_d \le d\hskip 2pt \frac{q-1}{q}\hskip 2pt q^{\ell+\frac{d}{\ell}}$ (as $\displaystyle \ell + \frac{d}{\ell}-\ell^\prime - \frac{d}{\ell^\prime} \ge 1$).

\subsubsection{2nd stage: uniqueness results} We will use Ritt's theorems to control the number of possible decompositions of a given polynomial.

\begin{proposition}
  \label{prop:uniq}
  Let $K$ be a field and $f\in K[x]$ be a polynomial of degree $d>0$. Assume the characteristic $p$ of $K$ does not divide $d$. Suppose we have two decompositions  $f = u \circ v = u^\prime
  \circ v^\prime$ of $f$ with 
  \begin{itemize}
  \item $u,v,u^\prime,v^\prime$ indecomposable,
  \item $\deg u = \deg u^\prime \ge 2$, $\deg v = \deg v^\prime \ge 2$,
  \item with $v,v^\prime$ monic with a zero constant term.
  \end{itemize}
\noindent
Then $u=u^\prime$ and $v=v^\prime$.
\end{proposition}

\begin{proof}
This follows from the first Ritt theorem \cite[\S 1.3 theorem 7]{Sc}
which more generally describes in which cases an equality $G_1 \circ
\cdots \circ G_r=  H_1\circ
\cdots \circ H_s$ with $G_i$, $H_j$ indecomposable of degree $>1$ may hold.
\end{proof}

As an immediate consequence, we obtain the case $d=p^2$ of theorem \ref{th:one} (a):
namely we have $\displaystyle D_{p^2} = D_{p,p} = \frac{q-1}{q}\hskip 2pt q^{2p}$.

\subsubsection{3rd stage: lower bounds for $D_{\frac{d}{\ell},\ell}$ and $D_{\ell, \frac{d}{\ell}}$}
\label{ssec:lower bounds}

\begin{lemma}
  \label{lem:low}
Assume $d$ is not a prime number.  Then we have
  $$D_{\ell,\frac{d}{\ell}} \ge \frac{q-1}{q}q^{\ell+\frac{d}{\ell}} \left( 1 -  \frac{d/\ell}{q^{\frac{d}{\ell}-\frac{d}{\ell^2}-\ell+1}}\right).$$
  And the same inequality holds for $D_{\ell,\frac{d}{\ell}}$ replaced
  by $D_{\frac{d}{\ell},\ell}$.
\end{lemma}

\begin{proof} We only give the proof for $D_{\frac{d}{\ell},\ell}$ as computations for $D_{\frac{d}{\ell},\ell}$ are the same.  In $\DD_{\ell,\frac{d}{\ell}}$ we will only count those polynomials $f$ which decompose as
  $f=u\circ v$ with $u$ and $v$ as in proposition \ref{prop:uniq}.  Then we obtain
  \begin{align*}
    D_{\ell,\frac{d}{\ell}}
    &\ge \frac{1}{q(q-1)} I_\ell \cdot I_{\frac{d}{\ell}} \\
    &\geq   \frac{1}{q(q-1)} N_\ell \hskip 2pt (N_{\frac{d}{\ell}}-D_{\frac{d}{\ell}}) \qquad \text{($D_\ell=0$ as $\ell$ is prime)} \\
    &= \frac{1}{q(q-1)}  (q-1)q^\ell \hskip 2pt \left( (q-1)q^\frac{d}{\ell} - D_{\frac{d}{\ell}} \right) \\
    &= \frac{q-1}{q}q^{\ell+\frac{d}{\ell}} \left( 1  - \frac{D_\frac{d}{\ell}}{(q-1)q^{\frac{d}{\ell}}}\right)\\
  \end{align*}
 \vskip -3mm
  
 If $d$ is the product of exactly $2$ primes then
  $D_\frac{d}{\ell}=0$ and
  \vskip -4mm
  $$ D_{\ell,\frac{d}{\ell}} \ge \frac{q-1}{q}q^{\ell+\frac{d}{\ell}}\leqno{\hbox{\rm (*)}}$$
which in this case is better than the announced result.  
 
 If $d$ is the product of at least $3$ primes, use the practical upper bound for $D_d$
 obtained in \S \ref{ssec:first stage} to write $D_{\frac d \ell} \le \frac{d}{\ell}\hskip 2pt \frac{q-1}{q}\hskip 2pt q^{\ell+\frac{d}{\ell^2}}$
and deduce

$$ D_{\ell,\frac{d}{\ell}} \ge 
\frac{q-1}{q}q^{\ell+\frac{d}{\ell}} \left( 1  - \frac{(d/\ell)\frac{q-1}{q}q^{\ell+\frac{d}{\ell^2}}}{(q-1)q^{\frac{d}{\ell}}}\right)= 
    \frac{q-1}{q}q^{\ell+\frac{d}{\ell}} \left( 1 -  \frac{d/\ell}{q^{\frac{d}{\ell}-\frac{d}{\ell^2}-\ell+1}}\right)$$
\end{proof}

\subsubsection{Estimating the multiple decompositions}
Next we write
$$D_d \ge \card(\DD_{\ell,\frac{d}{\ell}} \cup \DD_{\frac{d}{\ell},\ell}) =  D_{\ell,\frac{d}{\ell}} +D_{\frac{d}{\ell},\ell} - \card(\DD_{\ell,\frac{d}{\ell}} \cap \DD_{\frac{d}{\ell},\ell})$$

\noindent
In order to estimate $D_d$ we need to estimate the intersection.

\begin{lemma}
\label{lem:nonuni} We have 

$$\left\{\begin{matrix}
&\displaystyle \card(\DD_{\ell,\frac{d}{\ell}} \cap \DD_{\frac{d}{\ell},\ell}) 
  \le \frac{d}{\ell}\hskip 2pt q^{\frac{d}{\ell^2}+2\ell-1}\hfill \\
&\displaystyle D_d \ge 2\hskip 2pt \frac{q-1}{q}\hskip 2pt q^{\ell + \frac{d}{\ell}}\left(1- \frac{2d}{\ell}\frac{1}{q^{\frac{d}{\ell}-\frac{d}{\ell^2}-\ell+1}} \right)\\
\end{matrix}\right.$$
\end{lemma}

The lower bound for $D_d$ is the remaining inequality to be proved in theorem \ref{th:one} (b).
The more precise inequality (**) in the proof below will complete the proof of theorem 
\ref{th:one} (a) in the special case $d=pp^\prime$.

\begin{proof}[Proof of lemma \ref{lem:nonuni}] (a) If $\gcd(\ell,d/\ell)=1$ then 
 $\card( \DD_{\ell,\frac{d}{\ell}} \cap \DD_{\frac{d}{\ell},\ell}) \le q^5$.
 
Indeed let $f \in \DD_{\ell,\frac{d}{\ell}} \cap \DD_{\frac{d}{\ell},\ell}$
  and let $f = u\circ v$ be a decomposition  with $\deg u = \ell$ and
  $\deg v= d/\ell$.  We follow Ritt's second theorem (see
  \cite[\S 1.4, theorem 8]{Sc} and the notation there).  The
  hypotheses of that result are satisfied because the
  derivative $u^\prime$ of $u$ is non zero; otherwise $f^\prime=0$,
 and so $f\in \Ff_q[x^p]$ and the characteristic $p$ of $\Ff_q$ divides
 $d = \deg f$.  
  In first case of Ritt's second theorem we
  have $L_1\circ u = x^rP(x)^n$ and $v\circ L_2 = x^n$ (where $r \ge
  0,$\  $P\in \Ff_q[x]$\ and $L_1$, $L_2$ are linear functions).  In our situation we get $n =
  \frac{d}{\ell}$ and $\ell = r+\frac{d}{\ell}\deg P$. Then $\deg P = \frac{\ell^2-\ell r}{d} \le \frac{\ell^2}{d} < 1$ so $\deg P=0$, $L_1\circ u = x^\ell$ and $v\circ L_2 =
  x^{\frac{d}{\ell}}$.  Considering all possible linear functions yield
  at most $(q-1)^2q^2$ such decompositions.
 In second case of Ritt's second theorem we have $L_1\circ u =
  D_m(x,a^n)$ and $v\circ L_2 = D_n(x,a)$, $a\in \Ff_q$ (where $D_n(x,a)$
  denote Dickson's polynomials).  We here obtain $m = \ell$
  and $n = \frac{d}{\ell}$.  Considering all possible linear functions and
  all $a\in \Ff_q$ yield at most $(q-1)^2q^3$ such decompositions.
  Finally we obtain

\medskip 

\noindent
(*)  \centerline{$\displaystyle \card(\DD_{\ell,\frac{d}{\ell}}\cap \DD_{\frac{d}{\ell},\ell}) \le   (q-1)^2q^2 + (q-1)^2q^3 \le q^5$} 
\smallskip 

(b) If
  $\gcd(\ell,d/\ell)\neq 1$ then we have
  $\card( \DD_{\ell,\frac{d}{\ell}} \cap \DD_{\frac{d}{\ell},\ell}) 
  \le \frac{d}{\ell}\hskip 2pt q^{\frac{d}{\ell^2}+2\ell-1}$.

Indeed let $f \in \DD_{\ell,\frac{d}{\ell}} \cap \DD_{\frac{d}{\ell},\ell}$
  and let $f = u\circ v$ be a decomposition  with $\deg u = \ell$ and
  $\deg v= d/\ell$.  By Ritt's first theorem and because
  $\gcd(\ell,d/\ell)\neq 1$ either $u$ or $v$ is
  decomposable.  
  But as $\ell$ is a prime $\DD_\ell$ is empty and so 
  $v\in \DD_{\frac{d}{\ell}}$.  Thus we obtain
  \begin{align*}
\displaystyle
   \card(\DD_{\ell,\frac{d}{\ell}} \cap \DD_{\frac{d}{\ell},\ell})
    &\le N_\ell \hskip 2pt \frac{1}{q(q-1)} \hskip 2pt D_{\frac{d}{\ell}} \\
    &\le \frac{1}{q(q-1)} \hskip 2pt (q-1) \hskip 2pt q^\ell  \hskip 2pt  \frac{d}{\ell} \hskip 2pt q^{\ell + \frac{d}{\ell^2}} \quad \text{(end of \S \ref{ssec:first stage})} \\
    &\le \frac{d}{\ell} \hskip 3pt q^{\frac{d}{\ell^2}+2\ell-1} \\
  \end{align*}
\vskip -5mm  
  
The proof follows as for all $d>6$ we have $\displaystyle \frac{d}{\ell^2}+2\ell-1 \ge 5$.
\end{proof}


\end{document}